\documentclass[oneside,11pt]{amsart}
\usepackage{amsmath,amsfonts, latexsym,amssymb}
\newtheorem{theorem}{Theorem}[section]
\newtheorem{lemma}[theorem]{Lemma}

\newtheorem{proposition}[theorem]{Proposition}
\theoremstyle{definition}

\newtheorem{example}[theorem]{Example}
\newtheorem{remark}[theorem]{Remark}
\numberwithin{equation}{section}

 \makeatletter
\renewcommand*\subjclass[2][2000]{%
  \def\@subjclass{#2}%
  \@ifundefined{subjclassname@#1}{%
    \ClassWarning{\@classname}{Unknown edition (#1) of Mathematics
      Subject Classification; using '1991'.}%
  }{%
    \@xp\let\@xp\subjclassname\csname subjclassname@#1\endcsname
  }%
}
 \makeatother

\begin{document}

\title[Quasiconformal harmonic mappings]{On boundary correspondence of q.c. harmonic mappings between smooth Jordan domains}

\author{David Kalaj}
\address{University of Montenegro, faculty of natural sciences and mathematics,
Cetinjski put b.b. 81000, Podgorica, Montenegro}
\email{davidk@t-com.me}



\begin{abstract}
A quantitative version of an inequality obtained in
\cite[Theorem~2.1]{mathz} is given. More precisely, for normalized
$K$ quasiconformal harmonic mappings of the unit disk onto a Jordan
domain $\Omega\in C^{1,\mu}\ $ ($0<\mu\le 1$) we give an explicit
Lipschitz constant depending on the structure of $\Omega$ and on
$K$.  In addition we give a characterization of q.c. harmonic
mappings of the unit disk onto an arbitrary Jordan domain with
$C^{2,\alpha}$ boundary in terms of boundary function using the
Hilbert transformations. Moreover it is given a sharp explicit
quasiconformal constant in terms of the boundary function.
\end{abstract}
\maketitle

\section{Introduction and auxiliary results}
Let $A=\begin{pmatrix}
  a_{11} & a_{12} \\
  a_{21} & a_{22}
\end{pmatrix}.$ We will consider the
matrix norm:  $$|A|=\max\{|Az|: z\in \mathbf R^2, |z|=1\}\
$$ and the matrix function $$l(A)=\min\{|Az|: |z|=1\}.$$ Let
$w=u+iv:D\to G$, $D, G\subset\mathbf C$, have partial derivative at
$z\in D$. By $\nabla w(z)$ we denote the matrix $\begin{pmatrix}
  u_{x} & u_{y} \\
  v_{x} & v_{y}
\end{pmatrix}.$ For the matrix $\nabla
w$ we have \begin{equation}\label{opernorm}|\nabla w|=|w_z|+|w_{\bar
z}|\end{equation} and $$l(\nabla w)=||w_z|-|w_{\bar z}||,$$ where
$$w_z := \frac{1}{2}\left(w_x+\frac{1}{i}w_y\right)\text{ and } w_{\bar z} := \frac{1}{2}\left(w_x-\frac{1}{i}w_y\right).$$

A sense-preserving homeomorphism  $w\colon D\to G, $ where $D$ and
$G$ are subdomains of the complex plane $\mathbf C,$ is said to be
$K$-quasiconformal (K-q.c), $K\ge 1$, if $w$ is absolutely
continuous on a.e. horizontal and a.e. vertical line and
\begin{equation}\label{defqc} |\nabla w|\le K
l(\nabla w)\ \ \ \text{a.e. on $D$}.\end{equation} Notice that,
condition (\ref{defqc}) can be written as $$|w_{\bar z}|\le
k|w_z|\quad \text{a.e. on $D$ where $k=\frac{K-1}{K+1}$ i.e.
$K=\frac{1+k}{1-k}$ },$$ or in its equivalent form
\begin{equation}\label{december} \left|\frac{\partial w}{\partial r}\right|^2 +
\frac{1}{r^2}\left|\frac {\partial w}{\partial \varphi}\right|^2\le
\frac{1}{2}(K +\frac 1K)J_w\ \ (z=re^{i\varphi}),
\end{equation}
where $J_w $ is the Jacobian of $w$ (cf. \cite{Ahl}, pp. 23--24).
Finally the last is equivalent to:

$$\frac{1}{K}\le \left|\frac{\frac{r{\partial w}}{\partial
r}}{\frac{\partial w}{\partial \varphi}}\right|\le K.$$

This implies the inequality
\begin{equation}\label{december1}\frac{1}{r^2}(1+\frac{1}{K^2})\left|{\frac{\partial w}{\partial
\varphi}}\right|^2\le K J_w\ \ (z=re^{i\varphi}).\end{equation} A
function $w$ is called \emph{harmonic} in a region $D$ if it has
form  $w=u+iv$ where $u$ and $v$ are real-valued harmonic functions
in $D$. If $D$ is simply-connected, then there are two analytic
functions $g$ and $h$ defined on $D$ such that $w$ has the
representation
$$w=g+\overline h.$$

If $w$ is a harmonic univalent function, then by Lewy's theorem (see
\cite{l}), $w$ has a non-vanishing Jacobian and consequently,
according to the inverse mapping theorem, $w$ is a diffeomorphism.
If $k$ is an analytic function and $w$ is a harmonic function then
$w\circ k$ is harmonic. However $k\circ w$, in general is not
harmonic.

Let $$P(r,x-\varphi)=\frac{1-r^2}{2\pi (1-2r\cos(x-\varphi)+r^2)}$$
denote the Poisson kernel. Then every bounded harmonic function $w$
defined on the unit disc $\mathbf U:=\{z:|z|<1\}$ has the following
representation
\begin{equation}\label{e:POISSON}
w(z)=P[w_b](z)=\int_0^{2\pi}P(r,x-\varphi)w_b(e^{ix})dx,
\end{equation}
where $z=re^{i\varphi}$ and $w_b$ is a bounded integrable function
defined on the unit circle $S^1:=\{z:|z|=1\}$.

In this paper we continue to establish Lipschitz and co-Lipschitz
character of q.c. harmonic mappings between smooth domains. This
class contains conformal mappings. The conformal case is well-known
(\cite{Ko}, \cite{SW}, \cite{rmr}, \cite{G}, \cite{POW}) but it
seems only here we yield an explicit constant even for conformal
case.

The first result in the area of q.c. harmonic mappings was
established by O. Martio (\cite{Om}). Recently there are several
papers with deals with topic (\cite{MMM}-\cite{kalann},
\cite{pk}-\cite{MP}). See also \cite{wan} for the similar problem of
hyperbolic q.c. harmonic mappings of the unit disk.

It is worth to mention the following fact, q.c. harmonic mappings
share with conformal mappings the following property (a result of M.
Mateljevic and P. Pavlovic). This property do not satisfy hyperbolic
q.c. harmonic mappings of the unit disk onto itself.

\begin{proposition}\label{prope}
If $w = P[f]$ is a q.c. harmonic mapping of the unit disk onto a
Jordan domain $\Omega$ with rectifiable boundary, then $f$ is an
absolutely continuous function.
\end{proposition}

The proof can be found in \cite{MP}, \cite{pk} or \cite{Kalaj}. We
will use Proposition~\ref{prope} implicitly in our main
Theorems~\ref{main} and \ref{convex}.

 Some of the notations are taken
from \cite{mathz}. Let $\gamma\in C^{1,\mu}$, $0<\mu\le 1$, be a
Jordan curve such that the interior of $\gamma$ contains the origin
and let $g$ be the arc length parameterization of $\gamma$. Then
$|g'(s)|=1$.  Let
\begin{equation}K(s,t)=\text{Re}\,[\overline{(g(t)-g(s))}\cdot i
g'(s)]\end{equation} be a function defined on $[0,l]\times[0,l]$.
Denote by $K$ its periodic extension to $\mathbf R^2$
($K(s+nl,t+ml)=K(s,t)$, $m,n\in \mathbf Z$).

Since $K(s+nl,t+ml)=K(s,t)$, $m,n\in \mathbf Z$, it follows from
\cite[Lemma~1.1]{mathz} that
\begin{equation}\label{kernelar2}|K(s, t)|\le C_\gamma
d_\gamma(g(s),g(t))^{1+\mu},\end{equation} for
\begin{equation}\label{cgama}C_\gamma=\frac{1}{1+\mu}\sup_{s\neq
t}\frac{|g'(s)-g'(t)|}{(s-t)^\mu}\end{equation} and $d_\gamma$ is
the distance between $g(s)$ and $g(t)$ along the curve $\gamma$
i.e.
\begin{equation}\label{kernelar3}d_\gamma(g(s),g(t))=\min\{|s-t|,
(l-|s-t|)\}.\end{equation}

Using (\ref{kernelar2}) and following the same lines as in the
proof of \cite[Lemma~2.7]{mathz} we obtain the following lemma.

\begin{lemma}\label{jacabove}
Let $w=P[f](z)$ be a Lipschitz continuous  harmonic function between
the unit disk $\mathbf U$ and a Jordan domain $\Omega$, such that
$f$ is injective, and $\partial \Omega=f(S^1)\in C^{1,\mu}$. Then
for almost every $e^{i\varphi} \in S^1$ one has
\begin{equation}\label{jakobo}\limsup_{r \to 1-0} J_w(r e^{i\varphi}) \leq
C_\gamma|f'(\varphi)|\int_{-\pi}^\pi
\dfrac{d_\gamma(f(e^{i(\varphi+x)}),f(e^{i\varphi}))^{1+\mu}}{x^2}
dx,\end{equation} where $J_w$ denotes the Jacobian of $w$ at $z$,
and $f'(\varphi):=\frac{d}{d\varphi}f(e^{i\varphi})$.
\end{lemma}

A closed rectifiable Jordan curve $\gamma$ enjoys a $B-$ chord-arc
condition for some constant $B> 1$ if for all $z_1,z_2\in \gamma$
there holds the inequality
\begin{equation}\label{24march}
d_\gamma(z_1,z_2)\le B|z_1-z_2|.
\end{equation}
It is clear that if $\gamma\in C^{1,\alpha}$ then $\gamma$ enjoys a
chord-arc condition for some for some $B_\gamma>1$.

We will say that the q.c. mapping $f:\mathbf U\to \Omega$ is
normalized if $f(1)=w_0$, $f(e^{2\pi/3 i}) = w_1$ and $f(e^{-2\pi/3
i}) = w_2$, where ${w_0w_1}$, $w_1w_2$ and $w_2w_0$ are arcs of
$\gamma=\partial \Omega$ having the same length $|\gamma|/3$.

The following lemma is a quasiconformal version of
\cite[Lemma~1]{SW}. Moreover, here we give an explicit H\"older
constant $L_\gamma(K)$.
\begin{lemma}\label{newle}
Assume that $\gamma$ enjoys a chord-arc condition for some $B$. Then
for every $K-$ q.c. normalized mapping $f$ between the unit disk
$\mathbf U$ and the Jordan domain $\Omega=\mathrm{int}\gamma$ there
holds
$$|f(z_1)-f(z_2)|\le L_\gamma(K)|z_1-z_2|^\alpha$$ for $z_1,z_2\in
S^1$, $\alpha = \frac{1}{K(1+2B)^2}$ and $L_\gamma(K) =4
(1+2B)2^{\alpha}\sqrt{\frac{K|\Omega|}{\pi\log 2}}.$
\end{lemma}

\begin{proof}

For $a\in \mathbf C$ and $r>0$, $ D(a,r):=\{z:|z-a|<r\}$.
 It is clear that if $z_0\in S^1:=\partial \mathbf U$, then, because of normalization, $f(S^1\cap
\overline{D(z_0,1)})$ has common points with at most two of three
arcs ${w_0w_1}$, $w_1w_2$ and $w_2w_0$. (Here $w_0$, $w_1$, $w_2\in
\gamma$ divide $\gamma$ into three arcs with the same length such
that $f(1)=w_0$, $f(e^{2\pi i/3})=w_1$, $f(e^{4\pi i/3})=w_2$, and
$S^1\cap \overline{D(z_0,1)}$ do not intersect at least one of three
arcs defined by $1$, $e^{2\pi i/3}$ and $e^{4\pi i/3}$).

Let $k_\rho$ denotes the arc of the circle $|z-z_0|=\rho<1$ which
lies in $|z|\le 1$ and let $l_\rho=|f(k_\rho)|$.

Let $\gamma_\rho:= f(S^1\cap D(z_0,\rho))$ and let $|\gamma_\rho|$
be its length. Assume $w$ and $w'$ are the endpoints of
$\gamma_\rho$ i.e. of $f(k_\rho)$. Then $|\gamma_\rho| =
d_\gamma(w,w')$ or $|\gamma_\rho| = |\gamma| - d_\gamma(w,w')$. If
the first case hold, then since $\gamma$ enjoys the $B-$chord-arc
condition, it follows $|\gamma_\rho|\le B|w-w'|\le Bl_\rho$.
Consider now the last case. Let $\gamma_\rho' = \gamma\setminus
\gamma_\rho$. Then $\gamma_\rho'$ contains one of the arcs
${w_0w_1}$, $w_1w_2$ and $w_2w_0$. Thus $|\gamma_\rho|\le
2|\gamma_\rho'|$, and therefore
$$|\gamma_\rho|\le 2Bl_\rho.$$

On the other hand, by using \eqref{opernorm}, polar coordinates
and the Cauchy-Schwartz inequality, we have
\[\begin{split}l_\rho^2&=|f(k_\rho)|^2=\left(\int_{k_\rho}|f_zdz +f_{\bar z} d\bar z|\right)^2\\&\le \left(\int_{k_\rho}|\nabla f(z_0+\rho
e^{i\varphi})| \rho d\varphi\right)^2 \\&\le \int_{k_\rho}|\nabla
f(z_0+\rho e^{i\varphi})|^2 \rho d\varphi\cdot \int_{k_\rho} \rho
d\varphi.\end{split}\]

Since $l(k_\rho)\le 2\rho \pi/2$, for $r\le 1$, denoting
$\Delta_r=\mathbf U\cap D(r,z_0)$, we have

\begin{equation}\label{hen}\begin{split}\int_0^r\frac{l^2_\rho}{\rho}d\rho &\le
\int_0^r\int_{k\rho} |\nabla f(z_0+\rho e^{i\varphi})|^2 \rho
d\varphi d\rho\\& \le K \int_0^r\int_{k\rho} J_f(z_0+\rho
e^{i\varphi}) \rho d\varphi d\rho = \pi
A(r)K,\end{split}\end{equation} where $A(r)$ is the area of
$f(\Delta_r)$. Using the first part of the proof it follows that,
the length of boundary arc $\gamma_r$ of $f(\Delta_r)$ does not
exceed $2Bl_r$ which, according to the fact $\partial
f(\Delta_r)=\gamma_r \cup f(k_r)$, implies $|\partial f(\Delta_r)|
\le l_r+2Bl_r$. Therefore by the isoperimetric inequality $$ A(r)
\le \frac{|\partial f(\Delta_r)|^2}{4\pi}\le \frac{(l_r
+2Bl_r)^2}{4\pi} = l^2_r\frac{(1+2B)^2}{4\pi}.$$

Employing now \eqref{hen} we obtain
$$F(r):=\int_0^r\frac{l^2_\rho}{\rho}d\rho \le
Kl^2_r\frac{(1+2B)^2}{4}.$$ Observe that for $0<r\le 1$ there hold
the relation $rF'(r) = l^2_r$. Thus $$F(r)\le K rF'(r)
\frac{(1+2B)^2}{4}.$$ It follows that, for $$\alpha =
\frac{2}{K(1+2B)^2}$$ there holds $$\frac{d}{dr}\log (F(r)\cdot
r^{-2\alpha})\ge 0$$ i.e. the function $F(r)\cdot r^{-2\alpha}$ is
increasing. This yields

$$F(r)\le F(1)r^{2\alpha}\le K \frac{|\Omega|}{2\pi}
r^{2\alpha}.$$ Now there exists for every $r\le 1$ an
$r_1\in[r/\sqrt 2,r]$ such that
$$F(r)=\int_0^r\frac{l^2_\rho}{\rho}d\rho\ge \int_{r/\sqrt
2}^r\frac{l^2_\rho}{\rho}d\rho = l^2_{r_1}\log \sqrt 2. $$ Hence
$$l^2_{r_1}\le K\frac{|\Omega|}{\pi\log 2} r^{2\alpha}.$$ Thus if
$z$ is a point of $|z|\le 1$ with $|z-z_0|=r/\sqrt 2$, then
$$|f(z)-f(z_0)|\le (1+2B)l_r\le (1+2B)l_{r_1}.$$ Therefore $$
|f(z)-f(z_0)|\le H|z-z_0|^\alpha, $$ where $$H =
(1+2B)2^{\alpha}\sqrt{\frac{K|\Omega|}{\pi\log 2}}.$$

Thus we have for $z_1,z_2\in S^1$ the inequality
\begin{equation}\label{there}|f(z_1)-f(z_2)|\le
4H|z_1-z_2|^\alpha.\end{equation}
\end{proof}
\begin{remark}
By applying Lemma~\ref{newle}, and by using the M\"obius
transformations, it follows that, if $f$ is arbitrary conformal
mapping between the unit disk $\mathbf U$ and $\Omega$, where
$\Omega$ satisfies the conditions of Lemma~\ref{newle}, then
$|f(z_1)-f(z_2)|\le C(f,\gamma)K|z_1-z_2|^\alpha$ on $S^1$.
\end{remark}
\section{Quantitative bound for Lipschitz constant}
The aim of this section is to prove Theorem~\ref{main}. This is a
quantitative version of \cite[Theorem~2.1]{mathz}. Notice that,
the proof presented here is direct (it does not depend on
Kellogg's nor on Lindel\"of theorem on the theory of conformal
mappings (see \cite{G} for this topic)).
\begin{theorem}\label{main}
Let $w=P[f](z)$ be a harmonic normalized $K$ quasiconformal mapping
between the unit disk and the Jordan domain $\Omega$. If
$\gamma=\partial \Omega\in C^{1,\mu}$, then there exists a constant
$L=L(\gamma,K)$ (which satisfies the inequality \eqref{L} below)
such that
\begin{equation}\label{desired}|f'(\varphi)|\le L\text{ for almost every
$\varphi\in [0,2\pi]$},\end{equation}

and

\begin{equation}\label{equati}|w(z_1)-w(z_2)|\le KL|z_1-z_2|\,\,
\text{ for }z_1,z_2\in \mathbf U.\end{equation}
\end{theorem}

\begin{proof}
Assume first that $w=P[f]$ is Lipschitz and thus
$$\mathrm{ess}\sup_{0\le \theta\le 2\pi}|f'(\theta)|<\infty.$$ It
follows that
\begin{equation}\label{ediel}\frac{\partial w}{\partial \varphi}(z)
= P[f'](z).\end{equation} Therefore for $\varepsilon>0$ there exists
$\varphi$ such that
\begin{equation}\label{maxsup}\left |\frac{\partial w}{\partial
\varphi}(z)\right |\le \mathrm{ess}\sup_{0\le \theta\le 2\pi}
|f'(t)| =: L \le |f'(\varphi)| +\varepsilon.\end{equation}

According to \eqref{december1} and (\ref{jakobo}) we obtain:
$$(1+\frac{1}{K^2})|f'(\varphi)|^2\le \frac{\pi}{4}C_\gamma
K|f'(\varphi)|\int_{-\pi}^\pi
\dfrac{d_\gamma(f(e^{i(\varphi+x)}),f(e^{i\varphi}))^{1+\mu}}{x^2}
dx.$$ If $$C_2=\frac{\pi}{4}C_\gamma \frac{K^3}{1+K^2}$$ then
\begin{equation}\label{mda}\begin{split}L-\varepsilon&\le C_2\int_{-\pi}^\pi
\dfrac{d_\gamma(f(e^{i(\varphi+x)}),f(e^{i\varphi}))^{1+\mu}}{x^{1+\mu}}
\frac{dx}{x^{1-\mu}}\\&\le C_2\int_{-\pi}^\pi
\dfrac{d_\gamma(f(e^{i(\varphi+x)}),f(e^{i\varphi}))^{1+\mu-\beta}}{x^{1+\mu-\beta}
}L^\beta \frac{dx}{x^{1-\mu}}.\end{split}\end{equation} Thus
$$(L-\varepsilon)/L^\beta\le C_2\int_{-\pi}^\pi
\dfrac{d_\gamma(f(e^{i(\varphi+x)}),f(e^{i\varphi}))^{1+\mu-\beta}}{x^{1+\mu-\beta}
}\frac{dx}{x^{1-\mu}}.$$

Choose $\beta$: $0<\beta<1$ sufficiently close to 1 so that
$\sigma=(\alpha-1)(1+\mu-\beta)+\mu-1>-1$. For example $$\beta =
1-\frac{\mu\alpha}{2-\alpha},$$ and consequently $$\sigma=
\frac{\mu\alpha}{2-\alpha}-1.$$ From Lemma~\ref{newle} and
(\ref{24march}), letting $\varepsilon\to 0$, we get $$L^{1-\beta}\le
C_2\cdot (B_\gamma L_\gamma)^{1+\mu-\beta} \int_{-\pi}^\pi
x^{\sigma}dx=C_3,$$ and hence
\begin{equation}\label{eqmc}L\le
C_3^{1/(1-\beta)}=C_3^{\frac{2-\alpha}{\mu\alpha}}.\end{equation}

By \eqref{ediel} it follows that $$|zg'(z)-\overline{zh'(z)}|\le
L.$$

On the other hand, $$|\nabla w| = |g'|+|h'|$$ is subharmonic. This
follows that \begin{equation}\label{usetohil}|\nabla w(z)|\le
\max_{|z|=1}\{|g'(z)|+|h'(z)|\} \le
K\max_{|z|=1}\{|g'(z)|-|h'(z)|\} =KL.\end{equation}

This implies \eqref{equati}.

Using the previous case and making the same approach as in the
second part of theorem \cite[Theorem~2.1]{mathz} it follows that $w$
is a Lipschitz mapping. Now applying again the previous case we
obtain the desired conclusion.

\end{proof}

\begin{remark}\label{rem}
The previous proof yields the following estimate of a Lipschitz
constant $L$ for a normalized $K-$quasiconformal mapping between the
unit disk and a Jordan domain $\Omega$ bounded by a Jordan curve
$\gamma\in C^{1,\mu}$ satisfying a $B-$chord-arc condition.
\begin{equation}\label{L}L\le 4\pi\left(\frac{\pi}2
\frac{K^3}{1+K^2}
C_\gamma\frac{2-\alpha}{\mu\alpha}\right)^{\frac{2-\alpha}{\mu\alpha}}
\left\{4B (1+2B)\sqrt{\frac{K|\Omega|}{\pi\log
2}}\right\}^{\frac{2}{\alpha}},\end{equation} where $$\alpha =
\frac{1}{K(1+2B)^{2}}$$ and $C_\gamma$ is defined in
\eqref{cgama}. See \cite{MP}, \cite{pk}, \cite{MMM} and
\cite{trans} for more explicit (more precise) constants, in the
special case where $\gamma$ is the unit circle.

\end{remark}
\section{Boundary correspondence under q.c. harmonic mappings}

If $w=g+\overline h$ is a harmonic function then $$w_\varphi
=i(zg'(z)-\overline{zh'(z)})$$ is also harmonic. On the other hand
$$rw_r=zg'(z)+\overline{zh'(z)}.$$ Hence the function $rw_r$ is
the harmonic conjugate of $w_\varphi$ (this means that $w_\varphi
+ irw_r$ is analytic). The Hilbert transformation of $f'$ is
defined by the formula $$H(f')(\varphi)=-\frac 1\pi
\int_{0+}^\pi\frac{f'(\varphi+t)-f'(\varphi-t)}{2\tan(t/2)}\mathrm
dt$$ for a.e. $\varphi$ and $f'\in L^1(S^1)$. The facts concerning
the Hilbert transformation can be found in (\cite{ZY}, Chapter
VII).

There holds

\begin{equation}w_\varphi =P[f']\text{ and }rw_r
=P[H(f')],\end{equation} if $w_\varphi$ and $rw_r$ are bounded
harmonic.


The following theorem provides a necessary and a sufficient
condition for the harmonic extension of a homeomorphism from the
unit circle to a $C^{2,\mu}$ Jordan curve $\gamma$ to be a q.c
mapping, once we know that its image is $\Omega = \mathrm{int}\,
\gamma$. It is an extension of the corresponding result
\cite[Theorem~3.1]{mathz} from convex domains to arbitrarily smooth
domains.
\begin{theorem}\label{convex}
Let $f:S^1\to \gamma$ be an orientation preserving absolutely
continuous  homeomorphism of the unit circle onto the Jordan curve
$\gamma=\partial \Omega \in C^{2,\mu}$. If $P[f](\mathbf U)=\Omega$,
then $w=P[f]$ is a quasiconformal mapping if and only if
\begin{equation}\label{first}0<l(f):=\mathrm{ess\,inf\,} l(\nabla w(e^{i\varphi})),\end{equation}
\begin{equation}\label{second}||f'||_\infty:= \mathrm{ess\,sup\,}|f'(\varphi)|<\infty\end{equation} and
\begin{equation}\label{third}||H(f')||_\infty:=\mathrm{ess}\sup_{\varphi}|H(f')(\varphi)|<\infty.\end{equation}

If $f$ satisfies the conditions \eqref{first}, \eqref{second} and
\eqref{third}, then $w=P[f]$ is $K$ quasiconformal, where
\begin{equation}\label{KK}K:=\frac{\sqrt{||f'||^2_\infty +
||H(f')||^2_\infty-l(f)^2}}{l(f)}.\end{equation} The constant $K$ is
the best possible in the following sense, if $w$ is the identity or
it is a mapping close to the identity, then $K=1$ or $K$is close to
$1$ (respectively).

\end{theorem}
\begin{proof}
Under the above conditions the harmonic mapping $w$, by a result of
Kneser, is univalent (see for example \cite[p. 31]{dur}). Therefore
$w=g+\overline{h}$, where $g$ and $h$ are analytic and $J_u =
|g'|^2-|h'|^2>0$. This infers that the second dilatation $\mu =
{h'}/{g'}$ is well defined analytic function bounded by $1$.

\subsection{\bf The proof of necessity} Suppose $w=P[f]=g+\overline h$ is a
$K-$q.c. harmonic mapping that satisfies the conditions of the
theorem. By \cite[Theorem~2.1]{kalann}) we have

\begin{equation}\label{equqc} |\partial w(z)|-|\bar\partial
w(z)|\geq \frac{C(\Omega, K,a)}K>0,\, z\in \mathbf U.
\end{equation}
By \cite[Thoerem~2.1]{mathz} or Theorem~\ref{convex} we get
\begin{equation}\label{be}|f'(\varphi)|\le L\, a.e.\end{equation} and \begin{equation}\label{beo}\lim_{r\to 1}|\partial
w(re^{i\varphi})|-|\bar\partial w(re^{i\varphi})|= |\partial
w(e^{i\varphi})|-|\bar\partial w(e^{i\varphi})|\ \
a.e..\end{equation} Combining \eqref{be}, \eqref{beo} and
(\ref{equqc}) we get (\ref{first}) and (\ref{second}).

Next we prove (\ref{third}). Observe first that
\begin{equation*}
w_r=e^{i\varphi}w_z+e^{-i\varphi}w_{\overline z}. \end{equation*}
Thus \begin{equation}\label{derrad}|w_r|\le |\nabla
w|.\end{equation} By using \eqref{derrad} and \eqref{usetohil} it
follows that
\begin{equation}
|w_r(z)|\le K L.\end{equation} The last inequality implies that
there exist the radial limits of the harmonic conjugate $rw_r$
a.e. and
\begin{equation}\label{hilbert}
\lim_{r\to 1}rw_r(re^{i\varphi})=\lim_{r\to
1}w_r(re^{i\varphi})=H(f')(\varphi)\ a.e,
\end{equation} where $H(f')$ is the Hilbert transform of $f'$. Since $rw_r$ is a bounded harmonic function it
follows that $rw_r=P[H(f')]$, and therefore
$$||H(f')||_\infty=\mathrm{ess} \sup |H(f')(\varphi)|<\infty.$$
Thus we obtain (\ref{third}).

\subsection{\bf The proof of sufficiency}
We have to prove that under the conditions (\ref{first}),
 (\ref{second}) and
(\ref{third}) $w$ is quasiconformal. This means that we need to
prove the function
\begin{equation}K(z)=\frac{|w_z|+|w_{\bar z}|}{|w_z|-|w_{\bar z}|} = \frac{1 + |\mu|}{1-|\mu|} \end{equation} is bounded.

Since $\mu= \overline{w_{\bar z}}/w_z$ is an analytic function it
follows that $|\mu|$ is subharmonic. (Notice that, as $\phi(t) =
\frac{1+t}{1-t}$ is convex this yields that $K(z) = \phi(|\mu(z)|)$
is subharmonic).

It follows from (1.1) that $w_\varphi$ is equals the
Poisson-Stieltjes integral of $f'$: $$w_\varphi(re^{i\varphi}) =
\frac{1}{2\pi}\int_0^{2\pi}P(r,\varphi-t)df(t).$$ Hence, by
Fatou's theorem, the radial limits of $f_\varphi$ exist almost
everywhere and $\lim_{r\to 1-} f_\varphi(re^{i\varphi}) =
f_0'(\theta)$ a.e., where $f_0$ is the absolutely continuous part
of $f$.

As $rw_r$ is harmonic conjugate of $w_\varphi$, it turns out that if
$f$ is absolutely continuous, then $$\lim_{r\to
1-}f_r(re^{i\varphi}) = H(f')(\theta)\,\, (a.e.),$$

and

$$\lim_{r\to 1-} f_\varphi(re^{i\varphi}) = f'(\theta).$$

As $$|w_z|^2+|w_{\bar z}|^2 = \frac 12 \left( |w_r|^2+
\frac{|f_\varphi|^2}{r^2}\right)$$ it follows that
\begin{equation}\label{boundir}
\lim_{r\to 1-}|w_z|^2+|w_{\bar z}|^2\le\frac{1}{2}(||f'||^2_\infty
+ ||H(f')||^2_\infty).
\end{equation}

To continue we make use of \eqref{first}. From \eqref{boundir} and
\eqref{first} we obtain that
\begin{equation}\label{ldef}\mathrm{ess} \sup_{\varphi\in[0,2\pi)}\frac{|w_z(e^{i\varphi})|^2+|w_{\bar z}(e^{i\varphi})|^2}
{(|w_z(e^{i\varphi})|-|w_{\bar z}(e^{i\varphi})|)^2}\le
\frac{||f'||^2_\infty +
||H(f')||^2_\infty}{2l(f)^2}.\end{equation}

Hence

\begin{equation}\label{from}|w_z(e^{i\varphi})|^2+|w_{\bar z}(e^{i\varphi})|^2 \le S
(|w_z(e^{i\varphi})|-|w_{\bar z}(e^{i\varphi})|)^2 \ \ \
(a.e.),\end{equation}

where

\begin{equation}\label{S}
S:=\frac{||f'||^2_\infty + ||H(f')||^2_\infty}{2l(f)^2}.
\end{equation}
According to \eqref{ldef}, $S\ge 1$. Let $$\mu(e^{i\varphi}):=
\left|\frac{w_{\bar z}(e^{i\varphi})}{w_z(e^{i\varphi})}\right|.$$

As $w$ is a diffeomorphism, $|\mu(e^{i\varphi})|\le 1$. Then
\eqref{from} can be written as follows: $$
1+\mu^2(e^{i\varphi})\le S(1-\mu(e^{i\varphi}))^2, $$ i.e. $\mu =
\mu(e^{i\varphi})$ satisfies the inequality
\begin{equation}\label{qineq}
\mu^2(S-1) -2\mu S + S - 1=(S-1)(\mu-\mu_1)(\mu-\mu_2)\ge 0,
\end{equation}

where $$\mu_1 = \frac{S+\sqrt{2S-1}}{S-1}$$ and $$\mu_2  =
\frac{S-1}{S+\sqrt{2S-1}}.$$

From \eqref{qineq} it follows that $\mu(e^{i\varphi})\le \mu_2$ or
$\mu(e^{i\varphi}) \ge \mu_1$. But $\mu(e^{i\varphi}) \le 1$ and
therefore

\begin{equation}\label{mu}
\mu(e^{i\varphi}) \le \frac{S-1}{S+\sqrt{2S-1}} \ \ \ \ (a.e.).
\end{equation}

As $\mu(z)=|a(z)|$, where $a$ is an analytic function, it follows
that

$$\mu (z) \le k:=\mu_2, $$ for $z\in \mathbf U$.

This yields that $$K(z) \le K := \frac{1+k}{1-k}=
\frac{2S-1+\sqrt{2S-1}}{\sqrt{2S-1}+1}=\sqrt{2S-1} ,$$ i.e.

$$K(z) \le \frac{\sqrt{||f'||^2_\infty +
||H(f')||^2_\infty-l(f)^2}}{l(f)}$$ which means that $w$ is $K
=\frac{\sqrt{||f'||^2_\infty + ||H(f')||^2_\infty-l(f)^2}}{l(f)}$
quasiconfomal. The sharpness of the last results follows from the
fact that $K=1$ for $w$ being the identity.
\end{proof}
\subsection{Two examples}
The following example shows that, a $K$ (with $K$ arbitrary close
to $1$) q.c. harmonic selfmapping of the unit disk exists, having
non-smooth extension to the boundary, contrary to the conformal
case.

\begin{example} (\cite{kalajpub}). Let $$
\theta(\varphi)=\varphi\frac{1+{b}\sin(\log|\varphi|-\pi/4)}{1+b\sin(\log\pi-\pi/4)},
\ \varphi \in [-\pi,\pi],$$ where $0<b<\sqrt 2/2$, and let
$w(z)=P[f](z)=P[e^{i\theta(\varphi)}](z).$ Then $w$ is a
quasiconformal mapping of the unit disc onto itself such that
$f'(\varphi)$ does not exist for $\varphi=0$.  Using a similar
approach as in Theorem~\ref{convex} it can be shown that
$$K_w:=\sup_{|z|<1} \frac{|w_z|+|w_{\bar z}|}{|w_z|-|w_{\bar z}|}
\to 1$$ as $b\to 0$ and this means that, there exists a q.c.
harmonic mapping close enough to the identity, but its boundary
function is not differentiable at $1$. Details we will discus
elsewhere.
\end{example}
The next example shows that, the condition \eqref{first} of the
main theorem is important even for harmonic polynomials.

\begin{example}
Let $w$ be the harmonic polynomial defined in the unit disk by:
$$w(z) =  {z-1 - (z-1)^2}+\overline {z -1} = 3z-3-z^2+\bar z.$$

Then  $w$ is a univalent harmonic mapping of the unit disk onto
the domain bounded by the $C^{\infty}$ convex curve $\gamma =
\{(4\cos t - \cos(2t) - 3, \sin(2t) - 2\sin(t)), t\in [0,2\pi)\}$.
But $w_z(1)=w_{\bar z}(1) = 1$, and therefore $w$ is not
quasiconformal.
\end{example}



\end{document}